\documentclass[reqno, 12pt,a4letter]{amsart}
\usepackage{amsmath,amsxtra,amssymb,latexsym, amscd,amsthm}
\usepackage{graphicx,color}
\usepackage[utf8]{inputenc}
\usepackage[mathscr]{euscript}
\usepackage{mathrsfs}
\usepackage[english]{babel}
\usepackage{color}
\newtheorem {theorem}{Theorem}[section]

\newtheorem {remark}{Remark}[section]

\newcommand{\R}{\Bbb R}
\newcommand{\B}{\Bbb B}
\newcommand{\C}{\Bbb C}
\newcommand{\dist}{{\rm dist}}

\def\EES{{\accent"5E e}\kern-.5em\raise.8ex\hbox{\char'23 }}
\def\ow{o\kern-.42em\raise.82ex\hbox{
   \vrule width .12em height .0ex depth .075ex \kern-0.16em \char'56}\kern-.07em}
\def\OW{o\kern-.460em\raise1.36ex\hbox{
\vrule width .13em height .0ex depth .075ex \kern-0.16em
\char'56}\kern-.07em}
\def\DD{D\kern-.7em\raise0.4ex\hbox{\char '55}\kern.33em}

\title{Global mixed \L ojasiewicz inequalities and asymptotic critical values}

\author{S\~i-Ti\d{\^e}p \DD INH$^\dagger$}
\address{Institute of Mathematics, VAST, 18, Hoang Quoc Viet Road, Cau Giay District 10307, Hanoi, Vietnam}
\email{dstiep@math.ac.vn}

\author{Krzysztof KURDYKA$^\ddagger$}
\address{Laboratoire de Math\'ematiques (LAMA) UMR-5127 CNRS, B\^atiment Chablais, Campus Scientifique, 73376 Le Bourget-du-Lac cedex, France}
\email{Krzysztof.Kurdyka@univ-savoie.fr}

\author{Ti\EES n-S\OW n PH\d{A}M$^*$}
\address{Department of Mathematics, University of Dalat, 1 Phu Dong Thien Vuong, Dalat, Vietnam}
\email{sonpt@dlu.edu.vn}

\thanks{$^\dagger$This author is partially supported by Vietnam National Foundation for Science and Technology Development (NAFOSTED) grant 101.04-2017.12 and the Vietnam Academy of Science and Technology (VAST)}

\thanks{$^*$This author was partially supported by Vietnam National Foundation for Science and Technology Development (NAFOSTED) grant number 101.04-2016.05.}

\keywords{\L ojasiewicz inequalities, asymptotic critical values}

\subjclass{}

\date{ \today}

\begin{document}
\maketitle

\begin{abstract} In this paper, we prove a version of global \L ojasiewicz inequality for $C^1$ semialgebraic functions and relate its existence to the set of asymptotic critical values.
\end{abstract}

\section{Introduction}

Let $f \colon \mathbb{R}^n \rightarrow \mathbb{R}$ be a $C^1$ semialgebraic function (i.e., its graph is a semialgebraic set in $\Bbb R^{n+1}$). For $A\subset\Bbb R^n$ and $x\in\Bbb R^n$, let $\dist(x , A)$ be the Euclidean distance of $x$ to $A$. By convention,  $\dist(x,\emptyset)=1$. Denote by $E(f)$ the set of $t \in \mathbb{R}$ for which there are no positive constants $c, \alpha$ and $\beta$ such that the following global \L ojasiewicz inequality holds:
\begin{equation}\label{old}
|f(x) - t|^{\alpha} + |f(x) - t|^{\beta} \ge c\, \dist (x, \{f = t\}) \quad \textrm{ for all } \quad x \in \mathbb{R}^n.
\end{equation}

Let
\begin{eqnarray*}
K_\infty(f) := \{ y \in {\Bbb R}
& | & \textrm{there exists a sequence $x^k \in \mathbb{R}^n$ such that}\\
&& \|x^k\| \rightarrow +\infty, \ f(x^k) \rightarrow y, \quad \|x^k\| \| \nabla f(x^k) \| \to 0 \}
\end{eqnarray*}
the {\it set of asymptotic critical values of }$f$. For the  case where $f$ is a polynomial function, an algorithm for computing this set is given in \cite{Kurdyka2014}. Moreover let
\begin{eqnarray*}
\widetilde{K}_\infty(f) := \{ y \in {\Bbb R}
& | & \textrm{there exists a sequence $x^k \in \mathbb{R}^n$ such that}\\
&& \|x^k\| \rightarrow +\infty, \ f(x^k) \rightarrow y, \quad \|\nabla f(x^k) \| \to 0 \}.
\end{eqnarray*}

\begin{remark}{\rm
(i) The set $K_\infty(f)$ is finite, but in general $\widetilde{K}_\infty(f)$ is not finite. See, for example, \cite{Ha2008,Kurdyka2000,Parusinski1997}.

(ii) By \cite[Theorems 2 and 3]{Dinh2014}, it follows that $E(f)\subset \widetilde K_\infty(f).$

(iii) The set $E(f)$ may be infinite; for example let $f(x, y) := \frac{x}{y^2 + 1},$ then $E(f)=\R$ and so $E(f)\not=K_\infty(f)$.

(iv) Suppose that $f$ is a polynomial. If $n  = 2$, remark that $E(f) \subset \widetilde K_\infty(f)\subset\widetilde K_\infty(f_\C)$ and since $\widetilde K_\infty(f_\C)= K_\infty(f_\C)$ by \cite{Ha2008}, where $f_\C$ is the complexification of $f$, we deduce that $E(f)$ is finite. If $n>2$, it may happen that $E(f)$ is infinite (see, for instance, \cite[Example 1.11]{Parusinski1997}). 
}\end{remark}

In this paper, we propose a version of \L ojasiewicz inequality by changing slightly the left-hand side of \eqref{old} such that the new inequality still holds for all but a finite number of values $t$. The validity of the new inequality is also related to the set of asymptotic critical values. In fact, we will prove the following result. 

\begin{theorem}\label{MixedLojasiewicz} 
Let $f \colon \mathbb{R}^n \rightarrow \mathbb{R}$ be a $C^1$ semialgebraic function. Assume that $t\not\in K_\infty(f)$. Then there exist some constants $\alpha>0$ and $c > 0$ such that 
\begin{equation}\label{new}
|f(x) - t|^{\alpha} + \|x\| |f(x) - t| \ge c\, \dist(x, V_t) \quad \textrm{ for all } \quad x \in \mathbb{R}^n.
\end{equation}
\end{theorem}

\section{Proof of the main result}
Without loss of generality, we may suppose that $t=0$ and from now on, we write $V$ instead of $V_0$. 

First of all, assume that $V=\emptyset.$ In this situation, it holds that
$$\inf_{x \in \mathbb{R}^n } [f(x)]^2 > 0.$$
In fact, if it is not the case, then we can see that 
$$\lim_{\tau \to +\infty} \min_{\|x\|^2 = \tau^2} [f(x)]^2 = 0.$$
Consequently, there exists an analytic curve $(R, +\infty) \rightarrow \mathbb{R}^n \times \mathbb{R}, \tau \mapsto (\varphi(\tau), \mu(\tau)),$ such that
\begin{enumerate}
\item[(a)] $f(\varphi(\tau)) \nabla f(\varphi(\tau)) = \mu(\tau) \varphi(\tau);$
\item[(b)] $\|\varphi(\tau)\| = \tau;$ and
\item[(c)] $\lim_{\tau \to +\infty} f(\varphi(\tau))  = 0.$
\end{enumerate}
Note that $f(\varphi(\tau)) \ne 0$ for all $\tau \ge R,$ so we can define $\lambda(\tau) := \frac{\mu(\tau)}{f(\varphi(\tau))}.$
Furthermore, we can write
\begin{eqnarray*}
f(\varphi(\tau)) &=& c \tau^\nu + \textrm{ lower order terms in } \tau,
\end{eqnarray*}
where $c \ne 0$ and $\nu < 0.$ We have
\begin{eqnarray*}
\frac{d (f \circ \varphi)(\tau)}{d \tau}
&=& \left \langle  \nabla f(\varphi(\tau)), \frac{d \varphi(\tau)}{d \tau} \right \rangle \\
&=& \lambda (\tau) \left \langle  \varphi(\tau), \frac{d \varphi(\tau)}{d \tau} \right \rangle.
\end{eqnarray*}
(The second equality follows from  Condition (a).) Hence,
$$2 \frac{d (f \circ \varphi)(\tau)}{d \tau} = \lambda (\tau) \frac{d\|\varphi(\tau)\|^2}{d \tau} = 2 \lambda (\tau) \tau = 2 \frac{\mu(\tau)}{f(\varphi(\tau))} \tau.$$
This, together with Conditions (a) and (b), implies that
$$\left | \frac{d (f \circ \varphi)(\tau)}{d \tau} \right | = \| \nabla f(\varphi(\tau)) \|.$$ 
It follows that
$$\| \nabla f(\varphi(\tau)) \| \|\varphi(\tau)\| = c \nu \tau^\nu + \textrm{ higher order terms in } \tau.$$ 
This combined with Condition (c) and $\nu < 0$ yields $0 \in K_\infty(f),$ which contradicts to our assumption. Therefore, 
$$\inf_{x \in \mathbb{R}^n } [f(x)]^2 > 0.$$
Consequently, there exists $\delta > 0$ such that $|f|^{-1}(s) = \emptyset$ for $s\leq\delta$, then for all $x\in\R^n$, we have $|f(x)|\ge\delta$. We prove that for all $\alpha>0$, there exists $c=c(\alpha)>0$ such that \eqref{new} holds. Indeed, for $\|x\|\le 1,$ we have $|f(x)|^\alpha\ge\delta^\alpha=\delta^\alpha \dist(x,V)$ and for $\|x\|\ge 1,$ we have $\|x\| |f(x)|\ge\delta=\delta \, \dist(x,V).$ Hence \eqref{new} holds for $c=\min\{\delta^\alpha, \delta\}$.

Now we assume that $V \ne \emptyset.$ We list the following facts :
\begin{itemize}
\item [(d)] Since $0\not\in K_\infty(f),$ there exist $c_0>0,$ $\delta>0,$ and $R>0$ such that 
\begin{equation}\label{0}\|x\| \|\nabla f(x)\|\ge c_0 \ \text{ for } \ \|x\|\ge R \ \text{ and }\ |f(x)|\le\delta.
\end{equation} 
Without loss of generality, we may assume that $\delta<\frac{c_0}{3}$ and $R \ge \dist(0,V)$.

\item [(e)] By \cite{Lojasiewicz1959,Lojasiewicz1965}, there exist constants $\alpha > 0$ and $c_1 > 0$ such that 
\begin{equation}\label{local}
|f(x)|^\alpha\ge c_1\, \dist(x,V) \ \text{ for } \ \|x\|\le 2R.
\end{equation}

\item [(f)] For each $x \in \mathbb{R}^n$ such that $|f(x)| \ge \delta$ and $\|x\|\ge 2R$, we have 
\begin{eqnarray*}
\|x\| |f(x)| & \ge & \delta\|x\| \ = \ \frac{2\delta}{3}  \frac{3}{2}\|x\| \ \ge \ \frac{2\delta}{3}\big(\|x\|+R\big)\\
& \ge & \frac{2\delta}{3}\big(\|x\| + \dist(0,V)\big) \ \ge \ \frac{2\delta}{3}\, \dist(x,V).
\end{eqnarray*} 
\end{itemize}
Now we consider the remaining case where $\|x\|\ge 2R$ and $|f(x)|\le\delta$. Assume that we have proved:
\begin{equation}\label{near}
\|x\||f(x)|\ge \frac{2c_0}{5}\, \dist(x,V).
\end{equation}
Of course, this, together with \eqref{local}, completes the proof of Theorem~\ref{MixedLojasiewicz}.

So we are left with proving \eqref{near}. By contradiction, assume that there exists $x^0$ such that $\|x^0\|\ge 2R$, $|f(x^0)|\le\delta$ and 
\begin{equation}\label{2}
\|x^0\| |f(x^0)|< \frac{2c_0}{5}\, \dist(x^0,V).
\end{equation}
It is clear that $f(x^0)\not=0$ so we have $0=\min_{x\in\R^n}|f(x)|<|f(x^0)|.$ We consider two cases:

\subsubsection*{Case 1: $\dist(x^0,V)\le\frac{\|x^0\|}{2}$} \

By Ekeland variational principle (see \cite[Corollary 11]{Ekeland1979}) with the data $\epsilon := |f(x^0)|$ and $\lambda := \frac{\dist(x^0,V)}{2}$, there exists $y^0$ such that 
\begin{eqnarray}
\label{3} && |f(y^0)|\le |f(x^0)|\le\delta, \\
\label{4} && \|x^0-y^0\|\le\lambda=\frac{\dist(x^0,V)}{2}\le\frac{\|x^0\|}{4}, \\
\label{5} &&  |f(x)|+\frac{\epsilon}{\lambda}\|x-y^0\|\ge |f(y^0)| \ \text{ for all }\ x\in\R^n. 
\end{eqnarray}

From \eqref{4}, it follows that $ \|x^0-y^0\|\le\lambda < \dist(x^0, V),$ and so $y^0 \not\in V$ and $f(y^0) \ne 0.$ Without loss of generality, we may assume that $f(y^0)>0$, then $f(x)>0$ for all $x$ close enough from $y^0$. Now \eqref{5} implies that $y^0$ is a local minimum of $f(x)+\frac{\epsilon}{\lambda}\|x-y^0\|.$ Consequently $0\in\nabla f(y^0)+\frac{\epsilon}{\lambda}\B^n$, where $\B^n$ denotes the unit closed ball in $\R^n.$ Hence by \eqref{2}, we have
$$\|\nabla f(y^0)\|\le\frac{\epsilon}{\lambda}=\frac{2|f(x^0)|}{\dist(x^0,V)} <\frac{4c_0}{5\|x^0\|}. $$
Hence
$$\|x^0\| \|\nabla f(y^0)\| < \frac{4c_0}{5}. $$
By \eqref{4}, we have
$$\|y^0\|\le\|x^0\|+\|x^0-y^0\| \le \|x^0\|+\lambda\le\frac{5\|x^0\|}{4}.$$
Consequently 
\begin{equation}\label{6}
\|y^0\| \|\nabla f(y^0)\| < c_0. 
\end{equation}
Note that, by \eqref{4},
$$\|y^0\|\ge\|x^0\|-\|x^0-y^0\|\ge\|x^0\|-\lambda\ge\frac{3\|x^0\|}{4}>R$$
and $|f(y^0)|\le\delta$ by \eqref{3}. So \eqref{6} contradicts to \eqref{0}.

\subsubsection*{Case 2: $\dist(x^0,V)>\frac{\|x^0\|}{2}$}\

By Ekeland variational principle (see \cite{Ekeland1979}) with the data $\epsilon := |f(x^0)|$ and $\lambda := \frac{\|x^0\|}{2}$, there exists $y^0$ such that
\begin{eqnarray}
\label{7} && |f(y^0)|\le |f(x^0)|\le\delta, \\
\label{8} && \|x^0-y^0\|\le\lambda=\frac{\|x^0\|}{2}, \\
\label{9} && |f(x)|+\frac{\epsilon}{\lambda}\|x-y^0\|\ge |f(y^0)| \ \text{ for all }\ x\in\R^n. 
\end{eqnarray}
Similarly to Case~1, we have
$$\|\nabla f(y^0)\|\le\frac{\epsilon}{\lambda}=\frac{2|f(x^0)|}{\|x^0\|},$$
which implies that
$$\|x^0\| \|\nabla f(y^0)\| \le 2|f(x^0)| \le 2\delta. $$
By \eqref{8}, we have 
$$\|y^0\|\le\|x^0\|+\|x^0-y^0\| \le \|x^0\|+\lambda=\frac{3\|x^0\|}{2}.$$
Therefore
\begin{equation}\label{10}
\|y^0\| \|\nabla f(y^0)\| \le 3\delta<c_0. 
\end{equation}
Note that, by \eqref{8},
$$\|y^0\|\ge\|x^0\|-\|x^0-y^0\|\ge\|x^0\|-\lambda=\frac{\|x^0\|}{2}\ge R$$
and $|f(y^0)|\le\delta$ by \eqref{7}. So \eqref{10} contradicts to~\eqref{0}.


\section{Some remarks}

\begin{itemize}
\item [(i)] For the class of $C^0$ semialgebraic functions, by replacing the gradient norm $\|\nabla f\|$ by the nonsmooth slope ${\frak m}_f$ (see e.g., \cite{Mordukhovich2006,Rockafellar1998}), Theorem \ref{MixedLojasiewicz}  still holds with the same proof. Note that by a Sard theorem for tame set-valued mappings with closed graphs (\cite{Ioffe2007}), the set of asymptotic critical values of $f$ is still finite.

\item [(ii)] If $f$ is a polynomial of degree $d$ in $n$ variables, by \cite{Acunto2005}, the exponent $\alpha$ can be made explicit by $\alpha:=\frac{1}{{\mathscr R}(n,d)}$ where ${\mathscr R}(n,d):=d(3d-3)^{n-1}$ if $d>1$ and ${\mathscr R}(n,d):=1$ if $d=1.$

\item [(iii)] The converse of Theorem~\ref{MixedLojasiewicz} does not always hold, i.e., Inequality~\eqref{new} may hold even if $t\in K_\infty(f)$ as we see in the following example. Consider the Broughton polynomial (see \cite{Broughton})
$$f(x,y):=x(xy-1)=x^2y-x.$$
We have three cases:
\begin{itemize}
\item [(a)] $|x|\le 1$ and $|y|\le 1$. Then by (i), there exists a constant $c_1>0$ such that 
\begin{equation*}\label{local1}
|f(x,y)|^\frac{1}{18}\ge c_1\, \dist\big((x,y),V\big) .
\end{equation*}

\item [(b)] $|x|\ge 1$. We have $\|(x,y)\|\ge 1$ and $y=\frac{f(x,y)+x}{x^2}$, so 
$$\dist\big((x,y),V\big)<\Big|\frac{f(x,y)+x}{x^2}-\frac{1}{x}\Big|=\Big|\frac{f(x,y)}{x^2}\Big|\le\|(x,y)\| |f(x,y)|.$$

\item [(c)] $|y|\ge 1$. We have $\|(x,y)\|\ge 1$ and by solving the equation $f(x,y)=x^2y-x$ with $x$ as variable, we get $x=\frac{1\pm \sqrt{1+4yf(x,y)}}{2y}$, so 
\begin{eqnarray*}
\dist\big((x,y),V\big)&\le&\max\Big\{\Big|\frac{1+\sqrt{1+4yf(x,y)}}{2y}-\frac{1}{y}\Big|,\Big|\frac{1- \sqrt{1+4yf(x,y)}}{2y}\Big|\Big\}\\
&=&\Big|\frac{1- \sqrt{1+4yf(x,y)}}{2y}\Big|\\
&=&\Big|\frac{1- (1+4yf(x,y))}{2y\big(1+\sqrt{1+4yf(x,y)}\big)}\Big|\\
&=&\Big|\frac{2f(x,y)}{\big(1+\sqrt{1+4yf(x,y)}\big)}\Big|\le 2\|(x,y)\| |f(x,y)|.
\end{eqnarray*}
We have finally $$|f(x,y)|^\frac{1}{18}+\|(x,y)\| |f(x,y)|\ge\min\{c_1,\frac{1}{2}\}\, \dist\big((x,y),V\big).$$
\end{itemize}

\item [(iv)] The statement of Theorem~\ref{MixedLojasiewicz} does not always hold if we replace $K_\infty(f)$ by $B_\infty(f)$, where $B_\infty(f)$ is the set of bifurcation values of $f$. Indeed, let 
$$f(X)=f(x,y,z):=z\Big(x^4+(xy-1)^2\Big).$$
It is clear that $f$ is a trivial fibration over $\R$ so $B_\infty(f)=\emptyset.$ Consider the following parameterized curve $s\mapsto X(s)=\big(\frac{1}{s},s,s\big),\ s\gg 1$. We have $\nabla f(x,y,z)=\Big(z(4x^3+2y(xy-1),2xz(xy-1),x^4+(xy-1)^2)\Big)$, so
$$\|X(s)\| \|\nabla f(X(s))\|=\|\big(\frac{1}{s},s,s\big)\| \|\big(\frac{4}{s^2},0,\frac{1}{s^4}\big)\|\sim s.\frac{1}{s^2}=\frac{1}{s}\to 0  \ \ \text{as} \ \ s\to\infty.$$
Moreover $f(X(s))=\frac{1}{s^3}\to 0$, so $0\in K_\infty(f)$. On the other side, since
$$\|X(s)\| |f(X(s))|\sim s \frac{1}{s^3}=\frac{1}{s^2}\to 0 \ \ \text{and} \ \ \dist(X(s),V)=s\to\infty,$$
there are no constants $\alpha, c$ such that 
$$|f(X)|^{\alpha} + \|X\| |f(X)| \ge c\,  \dist(X, V) \quad \textrm{ for all } \quad X \in \mathbb{R}^3.$$

\item [(v)] We can not put an exponent $\beta<1$ on $\|x\|$ in Inequality~\eqref{new} as we see as follows. Let 
$$f(x,y):=\frac{x}{\sqrt{y^2+1}},$$
and $t:=2.$ We have $\nabla f(x,y)=\big(\frac{1}{\sqrt{y^2+1}},\frac{-xy}{(y^2+1)^\frac{3}{2}}\big)$. So $\|(x,y)\| \|\nabla f(x,y)\|\ge\frac{\sqrt{x^2+y^2}}{\sqrt{1+y^2}}.$ Hence $\|(x,y)\|  \|\nabla f(x,y)\|\to 0$ if and only if $(x,y)\to (0,0)$. Consequently $K_\infty(f)=\emptyset$. Consider the following parameterized curve $s\mapsto X(s)=\big(\sqrt{s^2+1},s\big),\ s\gg 1.$ It is clear that $X(s)\in f^{-1}(1)$ for all $s$. Let $\B(X(s),\frac{s}{4})$ be the closed ball of radius $\frac{s}{4}$ centered at $X(s)$ and let $B(X(s),\frac{s}{4}):=\{(x,y)\ : \ |x-X(s)|\le\frac{s}{4}, \ |y-Y(s)|\le\frac{s}{4}\}$. Set $g(s):=\max_{(x,y)\in\B(X(s),\frac{s}{4})}f(x,y)))$. Then
\begin{eqnarray*}
g(s)\le\max_{(x,y)\in B(X(s),\frac{s}{4})}f(x,y)&=&\frac{\sqrt{s^2+1}+\frac{s}{4}}{\sqrt{(s-\frac{s}{4})^2+1}}\\
&=&\frac{\sqrt{s^2+1}+\frac{s}{4}}{\sqrt{\frac{9}{16}s^2+1}}\\
&=&\frac{\sqrt{1+\frac{1}{s^2}}+\frac{1}{4}}{\frac{3}{4}\sqrt{1+\frac{16}{9s^2}}}\to \frac{5}{3}<t.\\
\end{eqnarray*}
Consequently, for $s$ big enough, the ball $\B(X(s),\frac{s}{4})$ does not intersect the fiber $f^{-1}(t)$. Hence $\frac{s}{4}\le\dist(X(s),V_t).$ So 
$$\|X(s)\| |f(X(s))-t|=\|X(s)\|=\sqrt{s^2+1+s^2}<2s\le 8 \, \dist(X(s),V_t).$$
Therefore Inequality~\eqref{new} does not hold any longer if we put an exponent $\beta<1$ on $\|x\|$.
\end{itemize}

\section*{Acknowledgments} This research was partially performed while the authors had been visiting at Vietnam Institute for Advanced Study in Mathematics (VIASM). The authors would like to thank the Institute for hospitality and support. We would like also to thank professor J\'er\^ome Bolte for useful discussions that led to this paper.

\end{document}